\documentclass[11pt]{article}\topmargin=0cm\hfuzz=10pt  \oddsidemargin=0.1cm \textheight 210mm\textwidth=16.3cm
\usepackage{amsmath,amsfonts}
\usepackage{algorithmic}
\usepackage{algorithm}
\usepackage{array}
\usepackage[caption=false,font=normalsize,labelfont=sf,textfont=sf]{subfig}
\usepackage{textcomp}
\usepackage{stfloats}
\usepackage{url}
\usepackage{verbatim}
\usepackage{cite}

\usepackage{graphicx}
\usepackage{epsfig}
\usepackage{epstopdf}
\AppendGraphicsExtensions{.tif}

\usepackage[numbers]{natbib}
\usepackage{epsfig}

\usepackage{amssymb}

\usepackage{amsthm}
\newtheorem{theorem}{Theorem}[section]
\newtheorem{proposition}[theorem]{Proposition}
\newtheorem{lemma}[theorem]{Lemma}
\newtheorem{corollary}[theorem]{Corollary}

\newtheorem{definition}[theorem]{Definition}

\newtheorem{remark}[theorem]{Remark}
\usepackage{amssymb}
\newcommand{\comm}[1]{}
\usepackage{color}
\newcommand{\gr}[1]{\color{green}#1\color{black}}
\newcommand{\bl}[1]{\color{blue}#1\color{black}}

\def\e{\varepsilon}

\def\defi{\stackrel{{\scriptscriptstyle \Delta}}{=}}
\def\defi{:=}

\def\a{\alpha}

\def\o{\omega}
\def\O{\Omega}

\def\F{{\cal F}}
\def\w{\widehat}
\def \Ind{{\,\rm Ind\,}}
\def \Ind{{\mathbb{I}}}

\def\esssup{\mathop{\rm ess\, sup}}
\def\const{{\rm const\,}}

\def\R{{\bf R}}

\def\b{\beta}

\def\C{{\bf C}}

\def\ww{\widetilde}
\def\X{{\cal X}}

\def\oo{\bar}

\def\p{\partial}

\def\V{{\cal V}}

\def\BB{{\cal C}}

\newcommand{\be}{\begin{equation}}
\newcommand{\ee}{\end{equation}}
\newcommand{\bd}{\begin{displaymath}}
\newcommand{\ed}{\end{displaymath}}
\newcommand{\ba}{\begin{array}{ll}}
\newcommand{\ea}{\end{array}}
\newcommand{\baa}{\begin{eqnarray}}
\newcommand{\eaa}{\end{eqnarray}}
\newcommand{\baaa}{\begin{eqnarray*}}
\newcommand{\eaaa}{\end{eqnarray*}}
\font\sm=cmr10

\def\oo{\bar}

\def\a{\alpha}

\def\sinc{{\rm sinc}}

\def\ZZ{{\mathbb{Z}}}

\def\g{\gamma}

\def\ee{\e}
\def\WO{\stackrel{p}{W_r^1}}
\def\WOd{\stackrel{p}{W_r^d}}

\def\BBO{{\cal H}(O_1)}
 \def\aa{{\rm{a}}}
 \def\aa{a}

\def\BB{{\cal B}(\rho)}
\def\CC{{\cal C}(\rho)}
\def\TP{{\frac{1}{2\pi}}}

\def\Lrho{\rho^{-1}L_\infty(\R)}
\def\Linvrho {\rho^{} L_1(\R)}
\def\ellinvrho{\rho^{}\ell_1}
\def\ellrho{\rho^{-1}\ell_\infty}

\def\bl{}
\def\gr{\comm}
\date{Submitted:  September 12, 2024. Revised: April 9, 2025 }

\begin{document}

\title{ Sampling Theorem and  explicit interpolation formula   for {non-decaying} unbounded signals}

\author{Nikolai Dokuchaev }



\maketitle

\begin{abstract}
The paper establishes an analog Whittaker-Shannon-Kotelnikov sampling theorem for unbounded non-decaying band-limited signals. 
An explicit interpolation formula is obtained for signals  sublinear growth with rate of growth less than 1/2. At any time, the
rate of decay for the $k$th coefficients of this formula is $\sim 1/k^2$.  In addition, the paper obtains a method for calculating  the coefficients of
the interpolation formula applicable to signals with arbitrarily high rate of polynomial growth. 
\end{abstract}

\section{Introduction}
Problems of recovery of signals from incomplete observations
 were studied intensively  in different settings. This includes 
 recovering signals from samples.  It is known that smooth signals can be recovered arbitrarily small error from 
 samples taken with sufficient frequency.  This lead to the interpolation problem. There are many different  approaches 
 to this problem for smooth enough signals: splines, B-splines,  or mean values of signal between sampling points; see, e.g., 
 \cite{AldRev,AngKantor,SomeAnn,Bar1,BardaKantor,BMspline,Bouz,ButRev1,CostaKantor,FeiMD,Han,Li,NU},
 and the references therein. In particular, B-spline  methods have been applied   for multi-dimension  signals featuring polynomial growth, including the case 
defined on multidimensional spaces \cite{AldRev,AngKantor,Han, Li,NU}.   

In this paper, we focus on interpolation formulae such as 
\baa x(t)=\sum_{k=-\infty}^\infty \aa_k(t)x(k),
  \label{Gen}\eaa  
  where  $t\in\R$, $x:\R\to \C$ is a signal, $\aa_k(t)$ are some coefficients.

 For the Schwartz space  of rapidly decaying functions, there is an  interpolation formula \cite{RadExp}. 
 The coefficients of this formula have to  be calculated via numerical integration of some special function, so it can be called  semi-explicit;
 see also \cite{RamExp}. This formula includes the value of the Fourier transforms at the sampling points. 

For the band-limited signals  vanishing on $\pm \infty$, there is an interpolation formula  with explicitly given coefficients:
the classical Nyquist–Shannon sampling theorem  establishes that a band-limited signal can be  recovered without error from a discrete sample taken with a sampling rate  that is at least twice the maximum frequency
of  the signal (the Nyquist critical rate).  In particular, a band-limited  signal $x(\cdot) \in L_2(\R)$ 
 with the  support of the Fourier transform $\int_{-\infty}^\infty e^{-i\o s}x(s)ds$   contained in the interval $[-\pi,\pi]$ 
  can be recovered from its sample $\{x(k)\}_{k=-\infty}^\infty$ as
  \baa x(t)=\sum_{k=-\infty}^\infty \frac{\sin(\pi(k- t))}{\pi (k-t)}x(k).
  \label{Shan}\eaa 
 This is celebrated Whittaker-Shannon-Kotelnikov  interpolation formula, also known as 
 Whittaker–Shannon interpolation formula, Shannon's interpolation formula, and Whittaker's interpolation formula. It can be observed that since 
 the coefficients of this  interpolation formula are decreasing as $\sim 1/k$, it covers only signals  such that $x(t)\to 0$ with a certain rate   
 as $|t|\to +\infty$.

There are many works devoted to generalization of the Sampling Theorem; see e.g. the reviews in 
\cite{Eldar,KPT,Tr,GenH,V01,Z}  and literature therein.  Unfortunately, the generalisations listed therein do not lead to any alternative explicit interpolation 
  formulae
for standard band-limited signals from $L_1(\R)$, and they not cover the setting with non-decaying signals.

As far as we know, there are only one another explicit interpolation formulae in the existing literature: it is Vaaler formula
\cite{VaaExp}; see some analysis in \cite{RamExp}. This formula is different from (\ref{Gen}) since it includes the value of derivatives at the sampling points. 
The coefficients of this formula decay as $\sim 1/k$, similarly to (\ref{Shan}).

 The present paper obtains interpolation formulae  similar to (\ref{Shan}) 
 but with faster decreasing coefficients. This  itself is a useful feature, since it  allows to reduce the error caused by 
 truncation that is unavoidable in calculations.  However, the main focus of this paper  is extension of the interpolation 
 formula on unbounded and non-vanishing  band-limited signals. 
 Let us explain why it is important. It can be observed that a unbounded non-vanishing  signal usually
 can be modified  to a signal from $L_1(\R)$ without any loss of information, 
 for example, by replacement  a signal $x(t)$ featuring polynomial growth  by $e^{-|t|}x(t)$.  However,  at least for the case of  signals from $L_{1}(\R)$, these damping transformations represent  the convolutions in the frequency domain,  with smoothing kernels. Unfortunately, a band-limited signal will be transformed
 into a non-band-limited one along the way.  For the general type  unbounded two-sided processes, one could expect a similar impact of the damping transformations on the spectrum.
Therefore, it was essential to  develop a special  approach for obtaining explicit
 interpolation formulae  applicable to unbounded non-decaying signals $x$ featuring polynomial growth  such that 
 \baaa\sup_{t\in\R} (1+ |t|)^{-\a}|x(t)|<+\infty\eaaa
 for some $\a\ge 0$.
 
The paper obtains an analog of Sampling Theorem and a  interpolation formula (\ref{Gen})  for band-limited non-decaying one-dimensional unbounded  signals featuring  polynomial
  growth and such that their spectrum is contained in a interval $I\subset (-\pi,\pi)$. We show that, for these signals, 
  where the series absolutely converge for each $t\in\R$. The coefficients  in (\ref{Gen}) are calculated as 
  \baa
\aa_k(t)= \frac{1}{2\pi}\int_{-\pi}^\pi E(t,\o)e^{-i\o k}d\o
\label{ak}\eaa
 for some regular enough functions $E(t,\o):\R\to\C$ such that $E(t,\o)=e^{i\o t}$ for $\o\in I$ (Theorem \ref{propE2}).
 The requirements for  the  regularity of these  function
 are  defined by the allowed rate of growth  for a class of admissible signals $x$.  An algorithm for calculations of $\aa_k(t)$ is provided
 (Lemma \ref{exd}  below).
  In particular, it is shown that, for any $\a\ge 0$,  
  the functions  $E$ can be selected such that $\aa_k(t)$  decay faster than $\sim 1/k^\nu$ for some  $\nu>1+\a$.

  Furthermore,  for the special  case of the signals featuring sublinear  growth with 
  $\a\in [0,1/2)$, we calculated  the coefficients $\aa_k(t)$ explicitly 
  (Theorem \ref{ThM}).  
  The $k$th coefficients for the corresponding new interpolation formula (\ref{a})  are decreasing as $\sim 1/k^2$. 
  A more general case   of the signals featuring   subquadratic  growth with   $\a\in [0,3/2)$, we found a semi-explicit
  representation for the coefficients $\aa_k(t)$ 
  via integrals of functions $e^{-i\o k}\o^l$ for $l=1,2,3$.
  (Theorem \ref{ThD2} and formula (\ref{aD2}) in Section \ref{SecM}). The corresponding $|\aa_k(t)|$ decay faster than  $\sim 1/k^{5/2}$. 
  It is also shown how to construct  the coefficients $\aa_k(t)$  for  the general case where $\a\ge 3/2$ such that  the corresponding values 
 $\aa_k(t)$.

Some numerical experiments are described in Section \ref{SecE}.

  \section{The main result: Sampling theorem and interpolation formula}
 \label{SecM}
 \subsection*{Some notations}
Let  $\R$ and $\C$, be the set of all real and complex numbers, respectively, and let $\ZZ$ be the set of all integers.

For $r\in[1,\infty)$, we denote by $L_r(\R)$  the standard space  of all functions  $x:\R\to \C$,  considered up to equivalency,  such that
$\|x\|_{L_r(\R)}\defi \left(\int_{-\infty}^{\infty}|x(t)|^r dt\right)^{1/r}<+\infty$.
We denote by $ L_{\infty}(\R) $ the standard space  of all functions  $x:\R\to \C$,  considered up to equivalency,  such that
$\|x\|_{ L_{\infty}(\R) }\defi \esssup_{t\in\R}|x(t)|<+\infty$.

We denote by  $C( \R )$  the standard linear  space of continuous functions $f:  \R \to\C$
with the uniform norm $\|f\|_C\defi \sup_t |f(t)|$.

 Let $\a>0$ be given.  Let $\rho:\R\to (0,1]$ be defined as 
\baa\label{rho}
\rho(t)=(1+|t|)^{-\a}.
\eaa

Let $\Lrho$ be the space of functions $f:\R\to \C$ such that
$\rho f \in L_\infty(\R)$, 
with the  norm $\|f\|_{\Lrho}\defi \|\rho f\|_{L_\infty(\R)}$.

Let $\Linvrho$ be the space of functions $f:\R\to \C$ such that
$\rho^{-1} f \in L_1(\R)$, 
with the  norm $\|f\|_{\Linvrho}\defi \|\rho^{-1} f\|_{L_1(\R)}$. 

\begin{definition}\label{defG}
For a Borel measurable set  $D\subset  \R $ with non-empty interior,   let  
$x\in \Lrho $ be such that  
$\int_{-\infty}^\infty x(t)y(t)dt=0$
for any  $y\in \Linvrho$ such that  $Y|_{ \R \setminus D}\equiv 0$, where $Y$ is the Fourier transform of 
$y$. In this case, we say that $D$ is a spectral gap of $x\in  \Lrho$.
 \end{definition} 

For $\O\in (0,+\infty)$, we denote by  $\V(\rho,\O)$  the set of all signals $x\in \Lrho$ with
 the spectral gap $\R\setminus (-\O,\O)$. We call these signals band-limited.

It was shown in \cite{D24p} that signals from $\V(\rho,\O)$ are continuous for any $\a\ge 0$.

We use the terms "spectral gap" and "band-limited" because, for  a signal  $x\in\Lrho$, Definition \ref{defG} 
means that $X(\o)=0$ for $\o\in D$, where $X$ the Fourier transform of $x$.
 The standard Fourier transform is not applicable for general type non-vanishing signals from $\Lrho$, however, we will use the terms "spectral gap" and "band-limited" for them as well. It is shown in Section \ref{ssecBL} below that this is still justified with respect to the spectral properties of these signal.

We consider below spectral representation and some applications for unbounded signals  from $\Lrho$.

\begin{proposition}\label{prop1} For any  $\a\ge 0$ and $\O\in (0,\pi)$, there exists a sequence of functions $\{\aa_k(\cdot)\}_{k\in\ZZ}\subset \Linvrho$
such that $\sup_{t\in\R}\sum_{k\in\ZZ}(1+|k|)^{\a}|\aa_k(t)|<+\infty$ and (\ref{Gen}) holds, i.e., $x(t)=\sum_{k\in\ZZ }\aa_k(t)x(k)$ for any $t\in\R$.
The corresponding series is   absolutely convergent  for each $t$.
\end{proposition}

In particular, it follows that, for any $t$,  the coefficients $|\aa_k(t)|$ are decreasing as $\sim 1/k^\nu$ for some  $\nu>1+\a$.

Lemma \ref{propE2} and  Lemma \ref{exd}  below provide a way to calculate  $\aa_k(t)$ for any $\a\in[0,+\infty)$

For the case where $\a\in [0,1/2)$,  it was possible to calculate   
 these functions explicitly, as is Theorem \ref{ThM} below.

 Let $\O\in (0,\pi)$ be given.
Let some  even integer number $N$ be selected such that \baaa\label{N}
N>\frac{\O}{\pi-\O}.
\eaaa
 It can be noted that we need to select larger $N$ for larger $\O$, i.e., $N\to +\infty$ as $\O\to\pi -0$.  

 For $t\in [N,N+1)$ and  $\tau=t-N\in[0,1)$,  let us select \baaa
g(t)=\frac{\pi \lfloor t\rfloor }{t}=\frac{\pi N}{N+\tau}.
\eaaa
Here  $\lfloor t\cdot=\min\{m\in\ZZ:\ t\ge m\}$.

It is easy to see that,  for any  $t\in [N,N+1)$, we have  $g(t)\ge \pi N(N+1)^{-1}>\O$ and $g(t)t=\pi N$.

Assume that the function $g(t)$ is extended periodically from $[N,N+1)$ to $\R$.  
This function is right-continuous. In addition,  $g(m)=\pi$ and $(t-m)g(t)=\pi N$   for any integer $m$ and any $t\in[N+m,N+m+1)$. 
\begin{theorem}\label{ThM}  If  $\a\in [0,1/2)$ in (\ref{rho}), then,  for  band-limited signal  $x\in  \V(\rho,\O)$, for any $t\in \R$
and $m=\lfloor t-N\rfloor$, i.e., $t\in[N+m,N+m+1)$, we have that (\ref{Gen}) holds, i.e., $x(t)=\sum_{k\in\ZZ }a_k(t)x(k)$ for $t\in\R$,
where $a_m(t)=1-\frac{g(t)}{\pi}$, and 
\comm{Atau[k]=N*sinc(gg*(ttD[k]-mmm))/((ttD[k]-tau1))}
\baa
a_k(t)=\frac{(t-m)\sin[ g(t)(k-m)]}{\pi (k-m)(k-t)},\qquad k\neq m. \label{a}
\eaa 
 The corresponding series is  absolutely convergent  for each $t$ and  is absolutely convergent  uniformly  on any bounded interval of $\R$.
 \end{theorem}

In particular, formula (\ref{a}) implies that $a_k(k)=1$, and that $a_k(l)=0$ for any integers $k$ and $l$ such that $k\neq l$.

It can be noted that, under the assumptions of Theorem  \ref{ThM},  we have that    
\baaa
a_k(t)
=\frac{(t-m)g(t)\sinc[ g(t)(k-m)]}{\pi(k-t)}
=\frac{N\sinc[ g(t)(k-m)]}{ k-t},\quad k\neq m. \label{alt2}
\eaaa
We used here that  $(t-m)g(t)=\pi N$.
  \begin{theorem}\label{ThD2}  If $\a\in [0,3/2)$ in (\ref{rho}), then, for  any band-limited signal  $x\in  \V(\rho,\O)$,
for any $t\in \R$
and $m=\lfloor t-N\rfloor$, we have that (\ref{Gen}) holds with
\baa
\aa_k(t)=\frac{\sin[ g(t)(k-m)]}{\pi(k-m-t)}+\frac{1}{\pi}\int_{g(t)}^\pi [P(t,\o)\cos(-k+m)-Q(t,\o)\sin(-k+m)]d\o, \label{aD2}
\eaa
were  \baaa
&&P(t,\o)\defi p_3(t)\o^3+p_2(t) \o^2+p_1(t)\o +p_0,\qquad Q(t,\o)\defi q_2(t) \o^2+q_1(t)\o +q_0,
\eaaa
and where 
\baaa
&&p_3(t)=1/3, \quad\qquad p_2(t)=-\frac{1}{2}(g(t)+\pi),\quad\qquad  p_1(t)=\pi g(t), \quad \\
&&p_0(t)=-\frac{1}{3}g(t)^3+\frac{1}{2}(g(t)+\pi)g(t)^2-\pi g(t)^2+1,
 \eaaa
\comm{ \baaa
&&R(t,\o)\defi \frac{1}{3}(\o^3-g(t)^3)-\frac{1}{2}(g(t)+\pi)(\o^2-g(t)^2)+\pi g(t)(\o-g(t))+1,
\eaaa}
and 
\baaa
&&q_0(t)=\frac{\pi-\g(t)}{\g(t)^2-\pi^2-2g(t)^2+2g(t)\pi},  \\&&q_1(t)=1-2q_0(t)g(t),\qquad \qquad q_2(t)=-q_0(t)\pi^2-p_1(t)\pi.
\eaaa
\comm{IBm<-function(w){-IB(-w)}
RB<-function(w){w^3/3-gg^3/3-(gg+pi)*(w^2-gg^2)/2+gg*pi*(w-gg)+1}; }
 The corresponding series is  absolutely convergent  for each $t\in \R$.
 \end{theorem}
 
 \gr{Tests: ER1<-function(w){cos(-(ttD[k]+mmm)*w+(tau1-mmm)*B1(w))};
aa1[k]=integrate(Vectorize(ER1),lower=gg,upper=pi,stop.on.error = TRUE)
aa2[k]=2*sin(gg*(ttD[k]-mmm))/(ttD[k]-tau1);

------

 Exzmples: aa2[k]=2*sin(gg*(ttD[k]-mmm))/(ttD[k]-tau1); cos(-(ttD[k]-mmm)*w+tau*B1(w))
aa1[k]=integrate(Vectorize(ER1),lower=gg,upper=pi,stop.on.error = TRUE)
tau1=mmm+M+tau}

 \section{Background for the proofs}
 \label{SecB}
 In this section, we outline  some results being used  in the proof for the main Theorem \ref{ThM}.
\subsection*{Some notations and  definitions for spaces of functions}
\comm{We denote by $\Ind$  the indicator function. }
We denote by  $\overline{z}$ the complex conjugation. We denote by $\ast$ the convolution 
\baaa
(h\ast x)(t)\defi\int_{-\infty}^\infty h(t-s) x(s)ds, \quad t\in\R.
\eaaa

For a Banach space $\X$, we denote by $\X^*$ its dual.

For $r\in[1,\infty)$, we denote by $\ell_r$ the set of all processes (signals) $x:\ZZ\to \C$, such that
$\|x\|_{\ell_r}\defi \left(\sum_{t=-\infty}^{\infty}|x(t)|^r\right)^{1/r}<+\infty$. We denote by $\ell_\infty$ the set of all processes (signals) $x:\ZZ\to \C$, such that
$\|x\|_{\ell_\infty}\defi \sup_{t\in\ZZ}|x(t)|<+\infty$.

For $k=1,2,...$, we denote by  ${\WOd}(-\pi,\pi)$ the Sobolev  space of functions $f: [-\pi,\pi]\to\C$
that belong to $L_2(-\pi,\pi)$ together with the distributional derivatives
up to the $k$th order, and such that $df^{l}(-\pi)/d\o=df^{k}(\pi)/d\o^k$ for $k=0,1,2,3,...,d-1$.

Let  $\BB$  be the space of continuous functions $f\in C( \R )$
with the  finite norm $\|f\|_{\BB}\defi \|\w f(\o)\|_{\Linvrho}$, where 
$\w f(\o)=\int_{-\infty}^\infty e^{-i\o s}f(s)ds$ \comm{i.e., 
$\ww f(\o)=2\pi \w f(\o)$,  where $\w f(\o)=\frac{1}{2\pi}\int_{-\infty}^\infty e^{-i\o s}f(s)ds$}
 is the Fourier transform of $f$.
By the choice of this norm, this is a separable Banach space that is isomorphic to $\Linvrho$, and the embedding 
$C(\R)\subset \BB$ is continuous.

In particular, the definition for $\BB$ implies that $Y\in\BB$ in Definition \ref{defG}.

Let $\ellinvrho$ be the space of sequences $x:\ZZ\to\C$ such that 
$\|x\|_{\ellinvrho}\defi \sum_{k\in\ZZ}\rho(k)^{-1}|x(k)|<+\infty$.

Let $\ellrho$ be the space of sequences $x:\ZZ\to\C$ such that 
$\|x\|_{\ellrho}\defi \sup_{k\in\ZZ}\rho(k)|x(k)|<+\infty$.

Let  $\CC$  be  the space of functions $f\in C([-\pi,\pi])$
with the  finite norm $\|f\|_{\CC}=\|\w f\|_{\ellinvrho}\defi\sum_{k\in\ZZ}\rho(k)^{-1}|\w f_k|$, where 
$\w f_k=\frac{1}{2\pi}\int_{-\pi}^\pi e^{-i k s}f(s)ds$ are the Fourier coefficients of $f$.
This is a separable Banach space that is isomorphic to $\ellinvrho$. 
In particular, if $\rho\equiv 1$, then functions from   $\CC$ have  absolutely convergent Fourier series on $[-\pi,\pi]$. 

\begin{lemma}\label{lemma1BB} 
Let  integer $d>\a+1/2$ and $r\in (1,2]$ be such that $r>(d-\a)^{-1}$. Then the following holds.
\begin{enumerate}
 \item[(i)]
 the embedding 
$ W_r^d (\R)\subset\BB$ is continuous, and 
\item[(ii)] the embedding 
$\WOd(-\pi,\pi)\subset\CC$ is continuous.
\comm{\item 
If $f\in \BB$ and $g\in\BB$, then $h=fg\in\BB$ and $\|h\|_{\BB}\le \|f\|_{\BB}\|g\|_{\BB}$.  }
\end{enumerate}
\end{lemma}
It can be noted that, since  $d-\a>1/2$ and hence $(d-\a)^{-1}<2$,  required $r$ exist.

\begin{lemma}\label{lemmaBBm}   Let $r$ be such as in Lemma \ref{lemma1BB}. 
Let $h\in W_r^1(\R)$ be such that $h(-\pi)=h(\pi)$, let $H$ be its Fourier transform, and let $m\in\ZZ$.  Then 
\  $H(t,\cdot) e^{i\cdot m}|_{[-\pi,\pi]}\in\CC$, and
$\rho(m)\|H(t,\cdot) e^{i\cdot m}|_{[-\pi,\pi]}\|_{\CC}\le \|h_{[-\pi,\pi]}\|_{\WO(-\pi,\pi)}$ for any $t$. 
\end{lemma}

  \begin{lemma}\label{lemma2BB} 
 Let  $h\in \Linvrho$, and let $H$ be its Fourier transform.  Then $H e^{i\cdot t}=(h\ast e^{i\o \cdot})(t)$, 
 $H e^{i\cdot t}\in\BB$, and $\rho(t)\|H e^{i\cdot t} \|_{\BB}\le \|h\|_{\Linvrho}$ for any $t\in\R$. In addition, 
the function  $H e^{i\cdot t}$ is continuous in the topology of $\BB$ with respect to $t\in\R$.
\end{lemma}

\subsection*{Spectral representation for signals from $\Lrho$} 
\label{ssecSpec}
We assume that each $X\in L_1( \R )$ represents an element of the dual space  $C( \R )^*$ such that
 $\langle X,f\rangle=\int_{-\infty}^\infty X(\o)f(\o)d\o$ for $f\in C( \R )$. We will use the same notation
  $\langle \cdot,\cdot\rangle$ for the extension of this bilinear form  
  on  $\BB^*\times \BB$.
 
 We will use the space $\BB$ and its dual $\BB^*$  to define formally a  spectral representation 
for $x\in \Lrho$  via $X\in\BB^*$
such that \baa
\langle X,\oo f\rangle=\int_{-\infty}^\infty x(t)\oo \varphi(t)dt
\eaa
 for any $f\in \BB$, where $\varphi\in \Linvrho$  is  
the Fourier transform for $f$.  We use notation $X=\F x$ for the corresponding operator $\F:\Lrho\to\BB^*$.

For the case where  $\rho\equiv 1$, this definition was used in Chapter VI in \cite{Katz}.  
In Chapter III in \cite{Kahane}, a similar definition was used for the Fourier transforms for 
pseudo-measures on $[-\pi,\pi]$ represented as elements of $\ell_\infty$; see also \cite{D24,D24p,Rei}.

  It can be noted that, for any $x\in\Lrho$, the corresponding $X=\F x$ is the weak* limit in $\BB^*$ as $t\to +\infty$ of the 
  sequence of functions 
$X_m(\o)\defi \int_{-m}^m e^{-i\o t} x(t)dt$ defined on $ \R $ (see \cite{D24p}).

\begin{proposition}\label{lemmahxy} Let $h\in \Linvrho$.
  For any  $x\in \Lrho$ and $X=\F x$, $H=\F h$, we have that 
 \baa\label{xX}
(h\ast x) (t) =\TP\langle X, H(\cdot)e^{i\cdot t} \rangle\quad \hbox{for}\quad t\in\R.
 \eaa
 \gr{NO NEED IN THIS PAPER: For any $X\in\BB^*$, there exists an unique up to equivalency process 
 $x\in \Lrho$ such that (\ref{xX}) holds for any $h\in \BB$ for all $t$. 
 For this process, we have that $\|x\|_{\Lrho}\le \|X\|_{\BB^*}$,
 and $\F x=X$.}
\end{proposition}

It can noted that, for $x\in L_1(\R)$, $\F x$ is the standard Fourier transform. In this case, $X\in C(\R)$.   Since  $He^{i\cdot t}\in \Linvrho \subset L_1(\R)$,
 the right hand part of  (\ref{xX})  is defined by standard way as $\TP\int_{-\infty}^\infty  X(\o)H(\o)e^{i\o t}d\o$. 

\subsection*{Spectral representation for band-limited signals from $\Lrho$} 
 \label{ssecBL}
   The following lemma connects Definition \ref{defG} with the spectral representation. 
\begin{lemma}\label{lemmaG} 
A signal  $x\in \Lrho$ has a spectral gap $D\subset \R$ if and only if  
$\langle  \F x,f\rangle =0$ for any $f\in \BB$ such  that $f|_{ \R \setminus D}\equiv 0$.
In this case, for any $f_1,f_2\in\BB$,
  if $f_1(\o)=f_2(\o)$ for all $\o\in[-\O,\O]$  then $\langle  \F x,f_1\rangle=\langle  \F x,f_2\rangle$.
\end{lemma}
 \par
\begin{proposition}\label{propE}
 Suppose that $\O_1>\O$ and  a function $E:\R\times\R\to \C$ be such that the following holds:
 \begin{itemize}
 \item[(i)]
 $E(t,\o)=e^{i \o t}$ for all $t\in\R$ and $\o\in [-\O_1,\O_1]$.
 \item[(ii)]
For any $t$,  $E(t,\cdot)\in\BB$ and $ \sup_{t\in\R}\rho(t)\|E(t,\cdot)\|_{\BB}<+\infty$.
\end{itemize} Then, for any 
 $x\in \V(\rho,\O)$ and $X=\F x$, we have that 
 \baa\label{xXE}
x (t) =\TP\langle X, E(t,\cdot) \rangle\quad \hbox{for}\quad t\in\R.
 \eaa
\end{proposition}
\section{Proofs} 
\subsection*{Proof of  background results}
\comm{ The proofs for the statements listed in Section \ref{ssecSpec}  can be found in \cite{D24p}.}
\par
{\em Proof of Lemma \ref{lemma1BB}}.
 The proof for Lemma \ref{lemma1BB}(i) can be found in \cite{D24p}.
 Let us prove  statement (ii).  
 Let $f\in \WOd(-\pi,\pi)$, and let $\w f_k$ be its Fourier coefficients. 
By the Hausdorff-Young inequality, we have for $r\in (1,2]$ that
\baaa
\left(\sum_{k\in\ZZ}(1+|k|)^{d q}| \w f_k|^q\right)^{1/q}\le C_r  \|f\|_{\WOd(-\pi,\pi)}
\eaaa 
for $q\defi (1-1/r)^{-1}$ and some $C_r>0$ that is independent on $f$. Further, we have that
$d^lf(\o)/d\o^l\in L_r(-\pi,\pi)$ for $l=0,1,...,d$, and
\baaa
\sum_{k\in\ZZ}\rho(k)^{-1}|\w f_k|&=& 
\sum_{k\in\ZZ}(1+|k|)^{\a}|\w f_k|= \sum_{k\in\ZZ}(1+|k|)^{\a-d}(1+|k|)^{d}|\w f_k|
\\
&\le& \left(\sum_{k\in\ZZ}
(1+|k|)^{(\a-d)r}\right)^{1/r}
\left(\sum_{k\in\ZZ}(1+|k|)^{d q}|\w f_k|^q\right)^{1/q} \le C_{d,r}\|f\|_{\WOd(-\pi,\pi)},
\eaaa
where
\baaa
C_{d,r}\defi  C_r \left(\sum_{k\in\ZZ}
\frac{ (1+|k|)^{\a r}}{(1+|k|)^{dr}}\right)^{1/r}<+\infty.
\eaaa
We used here,  that, since $d>\a+1/2$, we have that $d-\a>1/2$ and $1/(d-\a)<2$.

By the choice of $r$, we have that $r>1/(d-\a)$  and $r(d-\a)=dr-\a r>1$. 

This proves statement (ii).  This completes the proof of Lemma \ref{lemma1BB}. $\Box$
 \par
{\em Proof  of Lemma \ref{lemmaBBm}}. 
 Let $\w h(k)$ be the Fourier coefficients of  
$h_{[-\pi,\pi]}$. 
For $m\in\ZZ$, we have that
\baaa
H(\o) e^{i\o m}= e^{i\o m}\sum_{\nu\in\ZZ}\w h(k) e^{-i \o \nu}
=\sum_{\nu\in\ZZ}\w  h(\nu)  e^{-i \o (\nu-m)}=\sum_{\nu\in\ZZ}\w  h(\nu)  e^{i \o (m-\nu)}=\sum_{\nu\in\ZZ}\w   h(m+\nu)  e^{-i \o \nu}
\eaaa for $\o\in\R$.
Hence, for any $t\in\R$ and $C_{r,q}$ from the proof of Lemma \ref{lemma1BB}, 
 \baaa \|He^{i\o m}\|_{\CC}=
\sum_{\nu\in\ZZ}\rho(\nu)^{-1} |\w h(\nu+m)|=\sum_{\nu\in\ZZ}\rho(\nu-m)^{-1} |\w h(\nu)|
\le \rho(m)^{-1}\sum_{\nu\in\ZZ}\rho(\nu)^{-1} |\w h(\nu)|\\ =\rho(m)^{-1}\|\w h\|_{\ellinvrho}\le 
C_{r,q}\,\rho(m)^{-1}\|h\|_{\WO(-\pi,\pi)}.
\eaaa
We used here that $\rho(m)\rho(\nu)\le \rho(\nu-m)$.
This completes the proof of Lemma \ref{lemmaBBm}. $\Box$

{\em The proof  of Lemma \ref{lemma2BB}} can be found in \cite{D24p}.

{\em Proof of Proposition \ref{lemmahxy}.} 
Let $g(s)\defi \oo h(t-s)$,\comm{ $g(p)\defi \oo h(t-p)$, $g(t-q)\defi \oo h(q)$, q=t-p, p=t-q}
 i.e., $\oo h(s)=g(t-s)$, $s\in\R$. Let  $G=\F g$. By the definitions,
 \baaa
(h\ast x) (t)=\int_{-\infty}^\infty x(s)h(t-s)ds=
\int_{-\infty}^\infty x(s)\oo g(s)ds=
\TP\langle X, \oo G \rangle\quad \hbox{for}\quad t\in\R.
 \eaaa
We have that $G(\o)=\int_{-\infty}^\infty e^{-i\o s}g(t)ds$ and \baaa
\oo G(\o)=\int_{-\infty}^\infty e^{i\o s}\oo g(s)ds=\int_{-\infty}^\infty e^{i\o s}h(t-s)dt
=\int_{-\infty}^\infty e^{i\o (t-\nu)}h(\nu)d\nu=e^{i\o t}H(\o).
\eaaa
$\Box$

{\em The proof  of Lemma \ref{lemmaG}} follows immediately from the definitions.

{\em Proof of Proposition \ref{propE}}. 
Let $\BBO$ be the set of all $h\in\BB$ such that $h(t)=0$ if $|t|>\O_1-\O$.
Let $y\in \C$ be defined as 
 \baaa
y(t)=   \TP\langle X, E(t,\cdot) \rangle.  
\eaaa
 By Proposition  \ref{lemmahxy}  and by  the choice of $E$ and $X$, we have that \baaa
 (h\ast y)(t)=\TP \langle X, H E(t,\cdot )\rangle  \eaaa
for all $t$ and all $h\in\BBO$, $H=\F h$. 

 By Proposition  \ref{lemmahxy} and Lemma \ref{lemmaG}, we have that \baaa
 (h\ast x)(t)=(h\ast y)(t)  \eaaa
for all $t$ and all $h\in\BBO$. Hence $x=y$.
This completes the proof of Proposition \ref{propE}. $\Box$


\subsection{Proof of Proposition \ref{prop1}}


\begin{lemma}\label{lemmaE} 
 Suppose that there exists a function $E:\R\times\R\to \C$ satisfying the conditions of 
 Proposition \ref{propE} with 
 $\O_1\in (\O,\pi)$   such that, for any $t$, 
   $E(t,\cdot)|_{[-\pi,\pi]}\in\CC$ and $ \sup_{t\in\R}\rho(t)\|E(t,\cdot)|_{[-\pi,\pi]}\|_{\CC}<+\infty$.
Let
\baaa
\aa_k(t)\defi \frac{1}{2\pi}\int_{-\pi}^\pi E(t,\o)e^{-i\o k}d\o.
\eaaa
 Then $\{\aa_k(t)\}_{k\in\ZZ}\in \ellinvrho$ for all $t\in\R$,  $ \sup_{t\in\R}\rho(t)
 \|\{\aa_k(t)\}_{k\in\ZZ}\|_{\ellinvrho}<+\infty$, and any signal $x\in \V(\rho,\O)$ can be represented as
 (\ref{Gen})
 for $t\in [N,N+1)$, where the corresponding series is absolutely convergent for each $t$ as well as  in $\Lrho$.
 
 In addition, if  $E(t,\o)=\overline{E(t,-\o)}$  for all  $t$ and $\o$,  then $\aa_k(t)\in\R$.
\end{lemma}

\par
\gr{VERIFIED IN R: It can be noted that, under the assumption of Lemma \ref{lemmaE},
  we have that $\aa_k(k)=1$ and $\aa_k(l)=0$ for all $t\in\R$ and all integers $k$ and $l$, $l\neq k$.}

{\em Proof of Lemma \ref{lemmaE}}. 
Suppose that  $E(t,\o)$ and $\{\aa_k(t)\}_{k\in\ZZ}$ are  such as described.

By the assumptions,  $e^{i\o t}=E(t,\o)$ 
for $\o\in[-\O_1,\O_1]$.  By Proposition  \ref{lemmaE}, it follows that \baaa
 x(t)=\TP \langle X,E(\cdot, t)\rangle.       
\eaaa
In particular, \baaa
x(k)= \TP\langle X, E(k,\cdot ) \rangle,\quad k\in\ZZ.
\eaaa
By the choice of $\aa_k$, and by the  properties of Fourier series, we have that    \baaa
E(t,\o)=\sum_{k\in\ZZ} \aa_k(t) e^{i \o k},\quad t\in\R,\quad \o\in [-\pi,\pi],
\eaaa 
where the series are absolutely convergent for any $t$, $\o$.  
\bl{The condition that $E(t,\cdot)|_{[-\pi,\pi]}\in\CC$ implies that $\{\aa_k(t)\}_{k\in\ZZ}\in \ellinvrho$   and  that the sum
\baa \label{EE}
E(t,\cdot)=\sum_{k\in\ZZ} \aa_k(t) E(k,\cdot)
\eaa   converges in $\BB$ for any $t$.  The condition that $ \sup_{t\in\R}\rho(t)\|E(t,\cdot)|_{[-\pi,\pi]}\|_{\CC}<+\infty$
implies that $ \sup_{t\in\R}\rho(t)
 \|\{\aa_k(t)\}_{k\in\ZZ}\|_{\ellinvrho}<+\infty$.}

We have that    \baaa
e^{i\o t}=E(t,\o)=\sum_{k\in\ZZ} \aa_k(t) e^{i \o k}=\sum_{k\in\ZZ} \aa_k(t) E(k,\o),\quad t\in\R,\quad \o\in
[-\O_1,\O_1],
\eaaa 
where the series are absolutely convergent for any $t$, $\o$. 
 
By (\ref{EE}), it follows that
\baaa
x(t)=   \TP\langle X, \sum_{k\in\ZZ} \aa_k(t) E(k,\o)\rangle =\sum_{k\in\ZZ} \aa_k(t)\frac{1}{2\pi}\langle X,  E(k,\o)\rangle  
= \sum_{k\in\ZZ} \aa_k(t)x(k).\eaaa
This completes the proof of Lemma \ref{lemmaE}.
$\Box$

Let some $\O_1\in (\O,\pi)$ be selected such that \baa
\label{O1}g(t)\ge \O_1.
\label{O1}\eaa
Clearly, such $\O_1$  exists.   
\comm{It is easy to see that,  for any  $t$, we have that $g(t)\ge \pi N(N+1)^{-1}\ge \O_1$. }

\begin{lemma}\label{propE2} Let $r$ and $d$, be such as defined in Lemma \ref{lemma1BB}, given $\a>0$.
 Let  a function  $E_N:[N,N+1)\times\R\to \C$ be defined  
such that $ E_N(t,\cdot)\in W_r^d(\R)$ for any $t\in[N,N+1)$, and 
\baa
&& E_N(t,\o)=e^{i \o t},\quad \o\in [-g(t),g(t)],  \label{EW1}\\
 &&\sup_{t\in [N,N+1)}\|E_N(t,\cdot)|_{W_r^d(\R)}<+\infty, \label{EW2}
\\
&& \frac{d^l}{d\o^l}E_N(t,\pi)= \frac{d^l}{d\o^l} E_N(t,-\pi),\quad t\in\R,\quad l=0,1,...,d-1. \label{EW3}
\eaa  
Further, let the function  $E:\R\times\R\to \C$ be defined as
\baaa
E(t,\o)=E_N(t-m,\o) e^{i\o m}, \quad &\o\in \R,\quad t\in[N+m,N+m+1),\quad m\in\ZZ. 
\eaaa
Then the conditions of  Proposition \ref{lemmaE} hold for these  $E$.  
\end{lemma}

Clearly, if $E_N(t,\o)=\overline{E_N(t,-\o)}$  for all  $t$ and $\o$, then
 $E(t,\o)=\overline{E(t,-\o)}$  for all  $t$ and $\o$.

\par
{\em Proof of  Lemma \ref{propE2}}.  
It is easy to see that, by (\ref{EW1}), condition (i) of Proposition \ref{propE} is satisfied for selected  $E$. 
 We have that
$E_N(t,\cdot)\in\BB\subset W_r^1(\R)$ for any $t\in[N,N+1)$. By (\ref{EW2}), we have that
\baaa
\sup_{t\in[N,N+1)}\|E_N(t,\cdot)|_{[-\pi,\pi]}\|_{\BB}\le
\const\cdot \sup_{t\in[N,N+1)}
\|E_N(t,\cdot)\|_{W_r^1(\R)}<+\infty.
\eaaa 
 By  Lemma \ref{lemma2BB}(iii), we have for $u\in\BB$  that 
 \baaa
e^{i\cdot m}u\in \BB,\quad \|u\|_{\BB}\le\rho(m)^{-1}\|e^{i\cdot m}u\|_{\BB}.
\eaaa 
 
 By (\ref{EW3}), we have that 
$E_N(t,\cdot)|_{[-\pi,\pi]}\in \WO(-\pi,\pi)\subset \CC$ for any $t\in[N,N+1)$, and \baaa
\sup_{t\in[N,N+1)}\|E_N(t,\cdot)|_{[-\pi,\pi]}\|_{\CC}\le \const\cdot \sup_{t\in[N,N+1)}\|E_N(t,\cdot)|_{[-\pi,\pi]}\|_{\WO(-\pi,\pi)}<+\infty.
\eaaa 

By Lemma \ref{lemmaBBm}(iii), for any 
$v\in\CC$, we  have that \baaa  e^{i\cdot m}v\in \CC,\quad 
\|e^{i\cdot m}v\|_{\CC}\le \rho(m)^{-1}\|v\|_{\CC} .\eaaa
 Hence  conditions of Proposition \ref{propE} is satisfied  for $E(t,\o)$.  
 It follows    $E(t,\o)=\overline{E(t,-\o)}$  for all  $t$ and $\o$.
This completes the proof of Proposition \ref{propE2}.
$\Box$

\begin{lemma}\label{exd}
Let $d\in\ZZ$ be such that $d>\a+1/2$. Let $r\in (1/(d-a),2]$.
 A possible choice of $E_N(\cdot,t)|_{\{\o:\ |\o|\in[g(t),\pi]\}}$  satisfying conditions of Lemma \ref{propE2} is 
 the following.
 If $d$ is even, then  $E_N(t,\o)=P(t,\o)+i Q(t,\o)$, where $P(t,\o):[N,N+1)\times \R\to \R$ and $P(t,\o):[N,N+1)\times\R\to \R$ are  polynomials with respect 
 to $\o$ and such that
 all the following conditions hold.
\begin{enumerate}
\item $P(t,\cdot)$ is even, $Q(t,\cdot)$ is odd,   $P(t,g(t))=1$,  $Q(t,g(t))=0$,  $Q(t,\pi)=0$. 
 If $d\ge 2$, then  $P'_\o(t,\pi)=0$,  $Q'_\o(t,g(t))=t$, $Q'_\o(t,\pi)=0$.
\item If $d\ge 3$, then $P''_{\o \o}(t,g(t))=-t^2$,  $Q''_{\o \o}(t,g(t))=0$, $Q''_{\o \o}(t,\pi)=0$.
\item If $d\ge 2l$ for an integer $l\ge 2$ then 
\baaa
\frac{\p^{2l-1} P}{\p \o^{2l-1}}{ (t,\o)}\Bigl|_{\o=  g(t)}=0,\quad  \frac{\p^{2l-1} Q}{\p \o^{2l-1}}{ (t,\o)}\Bigl|_{\o=  g(t)}=
(-1)^{l-1}t^{2l-1}.\eaaa 
\item If $d\ge 2l+1$ for an integer $l\ge 2$ then 
\baaa
\frac{\p^{2l} P}{\p \o^{2l}}{ (t,\o)}\Bigl|_{\o=  g(t)}=(-1)^{l-1}t^{2l-1},\quad  \frac{\p^{2l} Q}{\p \o^{2l}}{ (t,\o)}\Bigl|_{\o=  g(t)}=0.
\eaaa 
\end{enumerate}
The corresponding coefficients $\aa_k(t)$ are real.
\end{lemma}
It can be noted that  the choice of function $\xi$ does not affect the values $\aa_k(t)$.  

\par
{\em Proof of Lemma \ref{exd}}. Clearly, we have that
\baaa
\frac{\p^k e^{it \o}}{\p \o^k}\Bigl|_{\o=\pm g(t)}=(it)^k e^{\pm it g(t)}
=(it)^k,\qquad k=1,2,...,
\eaaa
since $t g(t)=\pi N$ and $N$ is even integer, i.e.  $e^{\pm i g(t)t}=1$.
 It can be verified directly  that the  choice of $E_N$ ensures that
\baaa
&&E_N(t,\o)\Bigl|_{\o=\pm g(t)}=e^{it\o}\Bigl|_{\o=\pm g(t)}, \qquad E_N(t,\pi)=E_N(t, -\pi),\\
&&\frac{\p^k E_Nt,\o)}{\p \o^k}\Bigl|_{\o=\pm g(t)}=\frac{\p^k e^{it\o}}{\p \o^k}\Bigl|_{\o=\pm g(t)},\qquad
\frac{\p^k  E_N(t,\o) }{\p \o^k} \Bigl|_{\o=-\pi}
=\frac{\p^k  e^{it \mu(t,\o) }}{\p \o^k} \Bigl|_{\o=\pi},\\
&&k=0,1,2,...,d-1.\eaaa
 This ensures that $E_N$  satisfying conditions of Lemma \ref{propE2}. This completes the proof of  Lemma \ref{exd}. $\Box$

\begin{lemma}\label{lemmaShift} Let  $E$ and $\{\aa_k\}$ be selected as in Lemma \ref{lemmaE} and Lemma \ref{propE2}. 
 For any  $k,m\in\ZZ$, we have that 
 \baaa
\aa_k(t+m)=\aa_{k-m}(t).
 \eaaa 
 \end{lemma}
 \par
{\em Proof of Lemma \ref{lemmaShift}}. Let $N$  be defined as  in Lemma \ref{propE2},
and let $t=N+\tau$, 
 By the definitions, $E(t,\o)=e^{iM\o}E(\tau,t)$, where $M\in \ZZ$ is such that $\tau=t-M\in [N,N+1)$.
 Hence $E(t+m,\o)=e^{im\o}E(t,\o)$ and
 \baaa
\aa_k(t+m)=\frac{1}{2\pi}\int_{-\pi}^\pi e^{im\o}E(t,\o)e^{-i\o k}d\o=
\frac{1}{2\pi}\int_{-\pi}^\pi E(t,\o)e^{-i\o (k-m)}d\o=\aa_{k-m}(t). 
 \eaaa  This completes the proof of Lemma \ref{lemmaShift}.
$\Box$

\begin{corollary}\label{corrShift}  Let  $t$, $E$ and $\{\aa_k\}$ be such as  selected as in Proposition \ref{lemmaE} and Lemma \ref{propE2}.  
 Then for any $m\in\ZZ$, any signal $x\in \V(\rho,\O)$ can be represented, for $t\in[N+m,N+m+1)$, as
\baa
x(t)=\sum_{k\in\ZZ }\aa_{k-m}(t-m)x(k).
\eaa
The corresponding series is absolutely convergent for each $t$. 
\end{corollary}

 {\em Proof of Corollary  \ref{corrShift}}.  Let $\ww x(t)\defi x(t+m)$. It is easy to see that $\ww x\in\V(\rho,\O)$. Clearly, $x(t)=\ww x(t-m)$ for all $t$. By Proposition \ref{lemmaS1}, we have that
  \baaa
\ww x(s)=\sum_{k\in\ZZ } \aa_{k}(s)\ww x(k),\quad s\in [0,1).
\eaaa
Hence, for $t\in [m+N,m+N+1)$,
\baaa
x(t)=\ww x(t-m)=\sum_{k\in\ZZ } \aa_{k}(t-m)\ww x(k)
=\sum_{k\in\ZZ } \aa_{k}(t-m)x(k+m)=\sum_{k\in\ZZ } \aa_{d-m}(t-m)x(d).
\eaaa
This completes the proof of Corollary \ref{corrShift}.
 $\Box$

\begin{remark} Corollary \ref{corrShift} is
due to the particular choice of acceptable $E$. 
Possibly,  there exist  choices of $E$ such that  conditions required for Lemma \ref{lemmaE}  are satisfied but 
Corollary \ref{corrShift}  does not hold. 
\end{remark}

\par
{\em Proof of  Proposition \ref{prop1}} follows from  Lemma \ref{lemmaE} and  Corollary  \ref{corrShift}.  $\Box$

\subsection{Proof of   Theorem \ref{ThM}}

\gr{We will denote by $a_k(t)$ the corresponding coefficients $\aa_k(t)$ defined as in Proposition \ref{lemmaE}
with $E$ and $g$ defined by  Lemma \ref{propE2}.}

 \begin{lemma}\label{lemmaS1} Assume that  $\a\in [0,1/2)$ and that $E$ is selected 
 as in Lemma \ref{propE2}    with $d=1$, 
 $P(t,\o)\equiv  g(t)$, $Q(,\o)\equiv 0$  in Lemma \ref{exd}.
 Let $a_k(t)$ be the corresponding coefficients defined as in Proposition \ref{lemmaE}.
 Then,  for $t\in[N,N+1)$, we have that $\aa_0(t)=1-\frac{g(t)}{\pi}$ and
\baaa
a_k(t)=\frac{t\sin[ g(t)k]}{\pi k(k-t)},\quad k\neq 0.\label{a000}
\eaaa
 \end{lemma}
 It  can be noted that, since $g(t)t=\pi N$ for $t\in[N,N+1)$, we have that  \baaa
\aa_k(t)=\frac{g(t) t\, \sinc(g(t)k)}{\pi(k-t)}=\frac{N\sinc(g(t)k)}{k-t},\quad k\neq 0. \label{alt}
\eaaa

\par
{\em Proof of Lemma \ref{lemmaS1}}. 
Given that  $\a\in[0,1/2)$, we can select  $d=1$. 
 Clearly,   $\aa_k(N)=0$ for  $k\neq N$,  and 
$\aa_N(N)=1$, by 
the choice of $g(k)=\pi$ for $k\in\ZZ$. 

By the choice of $E_N$,  we have that
\baaa
\aa_k(t)=\frac{1}{2\pi}\int_{-\pi}^{\pi} e^{-i\o k} E_N(t,\o) d\o=\frac{1}{2\pi}(\a_k(t)+\b_k(t)), \label{ak1}
\eaaa
where  
 \baaa
&&\a_k(t)=\int_{-g(t)}^{g(t)} e^{-i\o k}  e^{i\o t} d\o,\qquad \label{alpha}
\\ &&\b_k(t)=\int_{-\pi}^{-g(t)} e^{-i\o k} e^{i g(t)t} d\o
+\int_{g(t)}^{\pi} e^{-i\o k}  e^{i g(t)t} d\o=\int_{-\pi}^{-g(t)} e^{-i\o k} d\o
+\int_{g(t)}^{\pi} e^{-i\o k} d\o.  \label{ab1}
\eaaa
We used here again  that $g(t)t=\pi N$.

 Consider the case where  $k\neq 0$. In this case,
 \baaa
\a_k(t)=\int_{-g(t)}^{g(t)} e^{-i\o k}  e^{i\o t} d\o=\frac{e^{ig(t)(t-k)}-e^{-ig(t)(t-k)} }{i(t-k)}=
\frac{e^{-ig(t)k}-e^{ig(t)k} }{i(t-k)}
=-\frac{2\sin(g(t)k)}{t-k},\label{ab}
\eaaa 
\baaa
\b_k(t)=
\frac{e^{ig(t)k}-e^{i\pi k} }{-ik}+\frac{e^{-i\pi k}-e^{-ig(t)k} }{-ik}=
\frac{e^{-ig(t)k} -e^{ig(t) k} }{ik}=-\frac{2\sin(g(t)k)}{k},
\eaaa
and
\baaa
a_k(t)=\frac{1}{2\pi}(\a_k(t)+\b_k(t)) =-\frac{1}{2\pi}\left(\frac{2\sin(g(t)k)}{t-k} + e^{i g(t)t}\frac{2\sin(g(t)k)}{k}\right).  
\label{a2}
\eaaa
Hence
\baa
a_k(t)=\frac{1}{2\pi}(\a_k(t)+\b_k(t)) =
-\frac{1}{\pi} \sin(g(t)k)\left(\frac{1}{t-k} +\frac{1}{k}\right) 
=
\frac{1}{\pi} \sin(g(t)k)\frac{t}{k(k-t)}.
 \label{a1}
\eaa
\par
For the case where  $k=0$, we have  
\baaa
\a_0(t)=\int_{-g(t)}^{g(t)}   e^{i\o t} d\o=-\frac{e^{ig(t)t}-e^{-ig(t)t} }{i t }=
-\frac{e^{iN\pi}-e^{-iN\pi} }{i t }. \label{ab2}
\eaaa
Since $t\ge N>0$, we have that $\a_0(t)=0$. Further, we have
  $\b_0(t)=2e^{i g(t)t}(\pi-g(t))=2(\pi-g(t))$. Hence  $a_0(t)=1-\frac{g(t)}{\pi}$.
 This completes the proof of Lemma \ref{lemmaS1}. $\Box$

{\em Proof of Theorem   \ref{ThM}}. The proof for the case where $t\notin [N,N+1)$
follows immediately from Lemma \ref{lemmaS1} and Corollary \ref{corrShift}.
 Convergence for any $t$ follows immediately from  
Proposition
\ref{propE2}, Lemma
\ref{lemmaS1}, and Corollary  \ref{corrShift}
since  $g(t-m)=g(t)$ for all $m\in\ZZ$. 

Let us prove uniform convergence on bounded intervals. It suffices to consider intervals 
$I=[-T,T]$, where $T>0$ can be arbitrarily large.

For $t\in\R$, $k\in\ZZ$, $n\in\ZZ$,   $n>2T+1$,  \comm{n>\max(2T+1,m)$} let $b_k(t)\defi a_k(t)x(k)$,  let $A(t,n)\defi \{k\in\ZZ:\ t -n\le k \le t +n\}$,
and let $x_{n}(t)\defi \sum_{k\in A(t,n)}
 b_k(t).$  We have  that
\baaa x(t)=x_{n}(t)+\sum_{k\notin A(t,n)} b_k(t).
\eaaa
For any $t\in\R$, we have that 
\baaa 
&&\sup_{t\in [-T,T]}
\sum_{k\notin A(t,n)} \left| b_k(t)\right|=\sup_{t\in [-T,T]}\sum_{k\notin A(t,n)}\left|\rho(k)^{-1}a_k(t)\rho(k)x(k)\right|
\\&&\le  \sup_{t\in [-T,T]}  \left(\sum_{k\notin A(t,n)}\rho(k)^{-1}|a_k(t)| \right)\sup_{k\in\ZZ}
\rho(k)   |x(k)|
 \le c_n(t)\|x\|_{\Lrho},
\eaaa
where  \baaa
c_n(t)\defi \sup_{t\in [-T,T]} \left(\sum_{k\notin A(t,n)}\rho(k)^{-1}|a_k(t)| \right) .
\eaaa
We have that $|k-t|>|k|-2T>0$ for $k\notin A(t,n)$ and
\baaa
|a_k(t)|\le \frac{N}{|k-m||k-t|}, \quad k\neq m,\quad  a_m(t)=1-\frac{g(t)}{\pi}\in (0,1).
\eaaa
Hence 
\baaa
c_n(t)\le \sup_{t\in [-T,T]} \left(\sum_{k\notin A(t,n),\ k\neq m}\rho(k)^{-1}
 \frac{N}{|k-m||k-t|}+\rho(m)^{-1}a_m(t)\right)\\
=N\sup_{t\in [-T,T]}\left(\sum_{k\notin A(t,n),\ k\neq m}
\frac{ (|1+|k|)^\a}{|k-m||k-t|} +\rho(m)^{-1}a_m(t)\right)\\
\le N\sum_{k\notin A(t,n),\ k\neq m}
\frac{ (|1+|k|)^\a}{|k-m|\,||k|-2T|}+\rho(m)^{-1}.
\eaaa
Hence  $\sup_{t\in[-T,T]}\| b_\cdot(t)\|_{\ell_1}<+\infty$ and  $x_n(t)\to x(t)$ as $n\to+\infty$ uniformly in $t\in[-T,T]$.
This completes the proof of Theorem   \ref{ThM}. 
$\Box$

{\em Proof of Theorem   \ref{ThD2}}. We  can select  $E_N$ as in Lemma \ref{exd}
with $d=2$.
 It is easy to verify directly that the assumptions of  Lemma \ref{exd} are satisfied with this choices.
Then the interpolation formula and its absolute convergence for any $t$ follows from  
Proposition
\ref{propE2}, Lemma
\ref{lemmaS1}, and Corollary  \ref{corrShift}.
This completes the proof of Theorem   \ref{ThD2}. $\Box$
\begin{remark}  The integrals in  Theorem   \ref{ThD2} can be calculated explicitly by parts  
as integrals of polynomials multiplied on $e^{-i\o k}$., and $\aa_k$ can be presented explicitly. 
The resulting  formulae will be tremendous, but could streamline calculations.
 In fact, or any positive $d\in \ZZ$,   one can, similarly to Theorems \ref{ThM}-\ref{ThD2}, provide 
formulas for $\aa_k$.
 \end{remark}
 \section{Some numerical experiments}
\label{SecE}
In some straightforward numerical experiments, 
we applied truncated interpolation with 

\begin{itemize} 
\item the classical interpolation formula (\ref{Shan}) covering decaying   signals,
\item our new  formula (\ref{a}) covering unbounded signals with a sublinear growth rate $\a<1/2$;
\item the formula with the coefficients calculated as in Theorem \ref{ThD2} 
covering unbounded signals with the growth rate $\a<3/2$.
\end{itemize}
 
 For these experiments, we used 
$\O=\frac{5\pi}{6}$.
The number  $N$ was selected as the smallest   even number such that $N>\O_1/(\pi-\O_1)$ with  $\O_1=(\O+\pi)/2$.
It can be noted that these choices define the coefficients of selected  interpolations formulae  uniquely.

We experimented with band-limited  signals from $L_2(\R)$ and unbounded signals featuring linear growth. 
 
\subsubsection*{Experiments with  band-limited  decaying signals from $L_2(\R)$} 
For band-limited  signals from $L_2(\R)$, we found that the results were
 indistinguishable for  all three interpolation  formulae.  

 First, let us show the results of our experiments  using an example of  a  band-limited signal $x(t)=\sinc(\O (t-1)+1/2)-2\sinc(\O(t+2)/\sqrt{2}-1)$.
 
 This type of signals has been used  for numerical examples in \cite{KPT}. 

We estimated $x(t)$ at  arbitrarily selected single points.  We calculated the error of interpolation
as $|x(t)-\w x(t)|$,
where $x(t)$ is the true value, $\w x(t)$ is the value obtained by the corresponding  truncated
interpolation formulae with summation over $k\in\{-L+\lfloor t\rfloor,-L+\lfloor t\rfloor+1,...,L+\lfloor t\rfloor\}$, for various integers $L>0$.

For  $t=-1.71$ and  $\O=5\pi/6$, we can  select   $N=6$ and $m=-8$. For this case, 
$x(t)= -1.8827743114725989937$. Following our approach, 
we select   $N=N(\O)=6$; this means that we have to select  $m=-8$. We obtained the following results.
\begin{itemize}
\item
For $L=50$,  the  error for the classical formula (\ref{Shan}) is 
 the  error for the classical formula (\ref{Shan}) is $7.3444679278278357515
 \times 10^{-6}$, 
the  error for our formulae  (\ref{Gen})-(\ref{a}) in Theorem \ref{ThM}  is
 $ 3.0762119742622218155\times 10^{-6}$, the  error for  our formulae  in Theorem \ref{ThD2} is 
 $2.3713406063219366615\times 10^{-6}$.
  \comm{
 [1] "A=1; xBL=function(t){1*sinc(Om*t-1*Om+1/2) -2*sinc(Om*(t+2)/sqrt(2)-1)}"
> T; Om; 
[1] 50
[1] 2.6179938779914944114
> tau1; tau; mmm;  true; #sigma
[1] -1.7099999999999999645
[1] 6.2900000000000000355
[1] -8
[1] -1.8827743114725989937
 > ErrorS
[1] 7.3444679278278357515e-06
> ErrorMyd1
[1] 3.0762119742622218155e-06
> ErrorMyd2
[1] -2.3713406063219366615e-06+0i
 }
 
For $L=100$,  the  error for the classical formula (\ref{Shan}) is
 $
 5.25190050757462501\times 10^{-5}$,  the  error for our formulae (\ref{Gen})-(\ref{a}) 
  in Theorem \ref{ThM}   is 
$2.4873476931475124729\times 10^{-6}$,
 the  error for  for our formulae  in  Theorem \ref{ThD2}  is $
5.0347767932557019321 \times 10^{-7}$.
\comm{
  ErrorS
[1] 5.25190050757462501e-05
> ErrorMyd1
[1] 2.4873476931475124729e-06
> ErrorMyd2
[1] 5.0347767932557019321e-07+0i}
\item
For $L=500$,  the  error for the classical formula (\ref{Shan}) is   $  2.2207117567063505703\times 10^{-7}$,
the  error for our formulae  (\ref{Gen})-(\ref{a}) in Theorem \ref{ThM}   is $4.3494727819393119717\times 10^{-9}$, 
 the  error for  for our formulae  in  Theorem \ref{ThD2}  is $
1.1088974183337541035\times 10^{-10}$.
 \end{itemize}

 \subsubsection*{Experiments with  band-limited unbounded  signals}
 In addition, we tested these interpolation formulae for some unbounded  signals. 
 In particular, we considered signal  
$x(t)=t\sin(\O (t-1)/1.0001+1/2)-2t\sin(\O(t+2)/\sqrt{2}-1)$ featuring linear growth.
This signals belong to $\Lrho$ with $\a=1$; this class is covered by   Theorem \ref{ThD2} but not  covered by the classical formula (\ref{Shan}) and
  by Theorem \ref{ThM}. 
    
  For $t=-1.71$, we have that 
  $x(t)=-1.0048884481864308604$,  and we have again that $N=6$, $m=-8$. We  found the following.
\begin{itemize}
\item
For $L=50$,  the  error for the classical formula (\ref{Shan}) is $\times 10^{-5}$,
the  error for our formulae  (\ref{Gen})-(\ref{a})
 in Theorem \ref{ThM}  is    $0.0045069457472912688445$, 
 the  error   for our formulae  in  Theorem \ref{ThD2} is $ 
 0.0053064993469169596807$.

 \item
For $L=100$,  the  error for the classical formula (\ref{Shan}) is $ 1.1277366571020146502$,  the  error for our formulae
 (\ref{Gen})-(\ref{a})  in Theorem \ref{ThM}   is  
$  0.05162016886519327219$,
 the  error for our interpolation calculated as  Theorem \ref{ThD2} is $0.011368950797937649178$.
\item
For $L=500$,  the  error for the classical formula (\ref{Shan}) is $0.16224049453604871829$,  
the  error for our formulae (\ref{Gen})-(\ref{a})  in Theorem \ref{ThM}  is $0.0013639041254578376794$,
 the  error for our formulae  in  Theorem \ref{ThD2}  is $ 8.1199227220052350162 \times10^{-5}$.
 \end{itemize}
 \comm{
 D1=AShen-aaa -red; D2=AShen-Atau -blue; D3=Atau-aaa - black 
 G=30
  plot(ttD[(T+1-G):(T+1+G)],D1[(T+1-G):(T+1+G)],
lwd = 1.0,cex=.5,type='p',col="red",xlab ='k',ylab='')# =expression(list(a[k](list(t)))))
 lines(ttD[(T+1-G):(T+1+G)],D2[(T+1-G):(T+1+G)],lwd = 1.5,cex=1.,type='p',col="light blue");
 lines(ttD[(T+1-G):(T+1+G)],D3[(T+1-G):(T+1+G)],lwd = 1.4,cex=.2,type='p',col="black");
 }

Let $\aa_k^{(0)}(t)$ be the coefficients for  the classical formula (\ref{Shan}), let $\aa_k^{(1)}(t)$ be the coefficients defined by 
   (\ref{Gen})-(\ref{a}) in Theorem \ref{ThM}, and 
 let   $\aa_k^{(2)}(t)$ be the coefficients calculated  as in  Theorem \ref{ThD2}.

Figure 1 shows  traces of the differences  
\baaa
\ww D_k(t)\defi  \aa_k^{(0)}(t)-\aa_k^{(2)}(t),\qquad  \oo D_k(t)\defi \aa_k^{(0)}(t)-a_k^{(1)}(t),\qquad 
D_k(t)\defi\aa_k^{(1)}(t)-\aa_k^{(2)}(t)
 \eaaa
for the coefficients $\aa_k(t)$ 
obtained by three different methods tested above. It can be seen 
that the differences between the coefficients are small. However, 
 truncated interpolation results are different as was shown above. 
Also, as expected, the asymptotic  behavior of the coefficients  $\aa_k(t)$ is different for these three methods. 
Figure 2 shows the values of the sequences
\baaa 
L_k\defi \log\left|\aa_k^{(1)}(t)|k|^{2.49}+1\right|,\quad  M_k\defi \log\left|\aa_k^{(2)}(t)|k|^{2.49}+1\right|,\quad \quad k=-500+\lfloor t\rfloor,...,500+\lfloor t\rfloor.
\eaaa 
The plot shows that the sequence $M_k$ is bounded, which is consistent 
with the results above implying  that $\aa_\cdot (t)\in \ellrho$ with $\a<3/2$, i.e., that  $\sum_{k\in\ZZ}(1+|k|)^{3/2}|\aa_k(t)|<+\infty$.
On the other hand,  the sequence $L_k$ is unbounded.  It can be noted that we used logarithmic  functions  here to make the trends more visible.

 \begin{figure}\label{fig1}\vspace{0cm}
\centerline{\psfig{figure=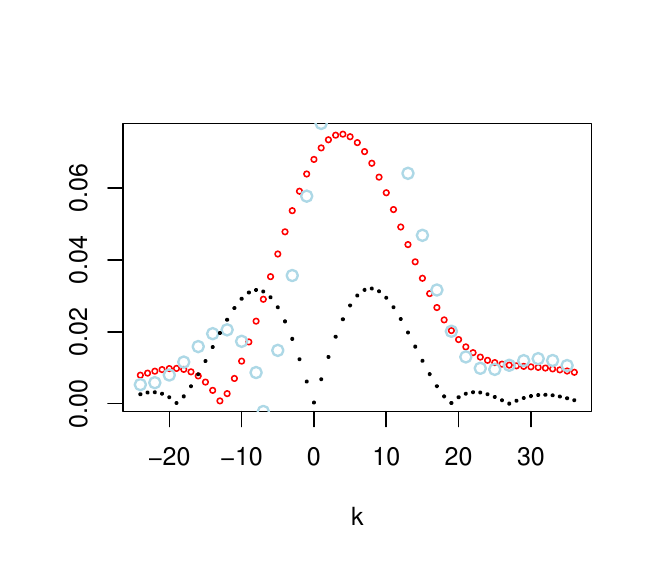,width=9cm,height=8.0cm}}
\vspace{-0.5cm}
\caption[]{\sm  The  differences $\ww D_k(t)$,  $\oo D_k(t)$, and $D_k$, between the  coefficients   $\aa_k(t)$ used in numerical experiments with $\O=5\pi/6$,  $t=-1.7$, and $N=6$, $m=-8$.  
Here big blue circles are for $\ww D_k(t)$,   medium red circles are for $\oo D_k(t)$, 
and black dots are for $D_k$.}
 \vspace{0cm}
 \end{figure}

 \begin{figure}\label{fig2}
\centerline{\psfig{figure=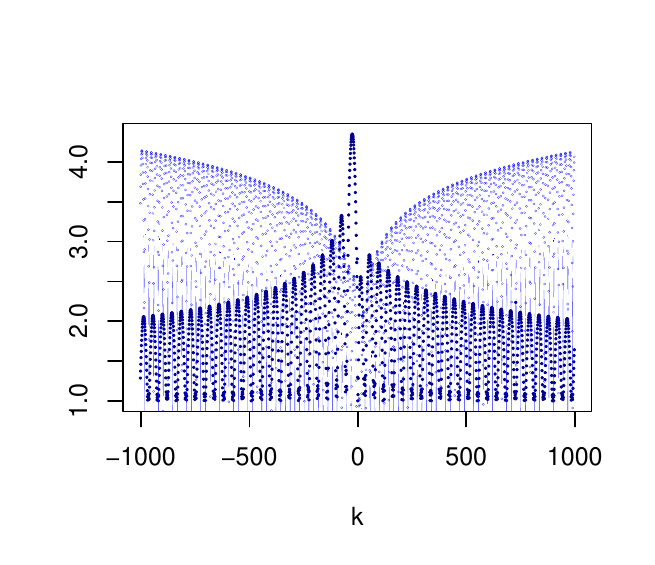,width=9cm,height=8.0cm}}
\vspace{-0.5cm}
\caption[]{\sm  The sequences $L_k(t)=\log\left(|k|^{2.49}|\aa_k^{(1)}(t)|+1\right)$ (light) and  $M_k=\log\left(|k|^{2.49}|\aa_k^{(2)}(t)|+1\right)$ (dark) 
for the coefficients  calculated  as   in  Theorem \ref{ThM}   and   Theorem \ref{ThD2}  respectively.}
 \vspace{0cm}
 \end{figure}

 \section*{Discussion and concluding remarks}
\begin{enumerate}
\item
The classical Whittaker-Shannon-Kotelnikov   interpolation formula (\ref{Shan})  is applicable
to bounded signals vanishing at $t\to \pm \infty$. Proposition \ref{prop1} covers unbounded 
 signals with  growth $\sim |t|^\a$ for any $\a\ge 0 $. 
 Theorem \ref{ThM} covers unbounded signals with sublinear growth $\sim |t|^\a$, $\a\in [0,1/2)$. 
  Theorem \ref{ThD2} covers unbounded signals with sub-quadratic  growth 
$\sim |t|^\a$, $\a\in [0,3/2)$.  
\item
For the case where $\a<1/2$, we have that  $|a_k(t)|\sim 1/k^2$ as $|k|\to +\infty$ and   $\sum_{k\in\ZZ}\rho(t)^{-1}|a_k(t)|<+\infty$ for any $t\in\R$.
This is why interpolation formula  (\ref{a}) is applicable to non-vanishing unbounded signals with moderate growth $\sim |t|^\a$. 
For these signals, 
the classical interpolation formula (\ref{Shan}) is not applicable, since its coefficients decay as $\sim 1/k$.
\item
The classical Whittaker-Shannon-Kotelnikov   interpolation formula (\ref{Shan})  allows spectrum bandwidth    
$[-\pi,\pi]$. On the other hand, Proposition \ref{prop1} and Theorem \ref{ThM} require that the spectrum bandwidth  of $x$ is  $[-\O,\O]$, 
for  $\O\in(0,\pi)$.  Therefore, the possibility to cover non-decaying unbounded signals is achieved via certain oversampling; 
this oversampling, however, can be arbitrarily small, since $\O$ can be arbitrarily close to $\pi$.
  
\item
 The condition  that $\O\in (0,\pi)$, and that the  sampling points are integers, can be removed, as usual,  by linear changes of the times scale, i.e., with the replacement of  the signal $x(t)$ by  signal $x(\mu t)$, with $\mu>0$. Clearly, less frequent sampling would require $\mu>1$,
and  selection of a larger $\O$  would require $\mu<1$.
\item
It can be seen that selection of $N$ used for the definition of $g$ in   interpolation formula  (\ref{Gen})-(\ref{a}) is non-unique.
\item
Different choices of $E(\cdot,\o)|_{|\o|>\O}$  
 lead to different  versions of the interpolation formula.
  
\item  We derived explicitly sampling coefficients $a_k(t)$ for the case where $\a\in (0,1/2)$.  
For the general  case where $\a>1/2$, one has to select more smooth functions  $E(\cdot,\o)$. 
Construction of  such functions is straightforward, as is shown in Lemma \ref{exd}(B). 
It is possible to obtain the corresponding coefficients $a_k(t)$ explicitely, similarly to the case where $\a<1/2$.
In particular, in the setting of Theorem \ref{ThD2} for $\a\in [0,3/2)$, this  would require calculation of  integrals of
$\o ^k e^{ -i\o k}$ for $k=1,2,3$. We leave it for the future research. Alternatively,  
these $\aa_k(t)$ this can be found via numerical integration; this is how it was done in our  numerical examples.
\item
It can be emphasised that the interpolation formulae  in Theorems \ref{ThM}-\ref{ThD2} are  exact;.
Therefore, for a vanishing signal $x\in L_2(\R)\cap \V(\rho,\O)$, all formulae from Theorem \ref{ThM}
and Theorem \ref{ThD2}  
give the same value as  (\ref{Shan}). Similarly, the value is the same  when 
applied for $x\in\V(\rho,\O)$  with $\aa_k$ calculated  for  $\O_0\in[\O,\pi)$ instead of $\O$. 
However, the values  for the corresponding finite truncated sums can be different.

\item  It can be shown that  Proposition \ref{lemmahxy} still holds for $\rho(t)=e^{-\a|t|}$ or $\rho(t)=e^{-\a t^2}$  with $\a>0$.
 However, it does not seem that analogs of Proposition
 \ref{prop1} or Theorems \ref{ThM}-\ref{ThD2} can be obtained in this case. This is because 
 all functions  $f\in\BB$ are analytical for these $\rho$, and hence if $f_{\R\setminus D}\equiv 0$ then $f\equiv 0$. 
This implies that  the exponents $e^{i\cdot t}$ do not feature analytical extension from $(-\O,\O)$ to $\R$  
 in the space $\Linvrho$ with these $\rho$.
Moreover, the concept of the spectral gaps  in Definition \ref{defG} 
is not applicable here;   this definition implies that any set $D\subset \R$   with non-empty interior, and such that $D\neq \R$, 
is a spectrum gap for any signal
 from $\Lrho$ with these $\rho$.

\end{enumerate}


\begin{thebibliography}{100}
\bibitem{AldRev}
Aldroubi, A., Krishtal, I, and Weber, E. (2015).  Finite dimensional dynamical sampling: an overview.
In: Excursions in Harmonic Analysis, Volume 4, pages 231-244. Springer.
\bibitem{AngKantor} 
 Angeloni, L., Cetin, N., Costarelli, D., Sambucini, A., and Vinti G. (2021) Multivariate sampling
Kantorovich operators: quantitative estimates in Orlicz spaces. Constructive Mathematical
Analysis, 4(2):229-241.
\bibitem{SomeAnn}
Annaby M., and P. Butzer, P. Sampling in Paley-Wiener spaces associated with fractional
integro-differential operators. J. Comput. Appl. Math., 171(1-2):39-57, 2004.
\bibitem{Bar1}
Bardaro, C.,  Butzer, P., and  Mantellini, I.(2014) The exponential sampling theorem of signal analysis
and the reproducing kernel formula in the Mellin transform setting. Sampl. Theory Signal
Image Process., 13(1):35-66.
\bibitem{BardaKantor} 
Bardaro, C.,Vinti, G., Butzer,  P.L., and Stens, R.L. (2007). Kantorovich-type generalized sampling
series in the setting of Orlicz spaces. Sampl. Theory Signal Image Process., 6(1):29, 2007.
\bibitem{BMspline} 
Bhandari, A., and P. Marziliano, P.  Fractional Delay Filters Based on Generalized Cardinal Exponential Splines.
Signal Processing Letters, IEEE 17(3), 225 - 228.
\bibitem{Bouz}
     Bouzeffour, F. A (2007). Whittaker-Shannon-Kotel’nikov sampling
theorem related to the Askey-Wilson functions. 
  J. Nonlinear Mathematical Physics 14 (3), 375–388 .
  \comm{A Whittaker-Shannon-Kotel’nikov sampling theorem related to the Askey-Wilson
functions is proved.
  There is an explicit interpolation formula for so-called  q-band-limited .
  To derive the sampling theorem associated to the Askey-Wilson function, we start by
defining a notion of the}

\bibitem{ButRev1}
Butzer, P.L., Stens, R.L. and Hauss, M. (1991). The sampling theorem and its unique role in various
branches of mathematics. Mitt. Math. Ges. Hamb., 12(3): pp. 523-547.
\bibitem{ButRev2}
 Butzer, P.. The sampling theorem of signal processing.  (2002). In Trends in industrial and applied
mathematics. Appl. Optim., 72,  p.  23-61. Kluwer Acad. Publ.,
Dordrecht.



\bibitem{CostaKantor} Costarelli, D.,  Minotti, A.M., Vinti. G. (2017).
Approximation of discontinuous signals by sampling Kantorovich series.  J. Math.Anal.Appl. 450, pp.1083–1103.
\bibitem{D24} Dokuchaev, N. (2024). 
Spectral representation of two-sided signals from  $\ell_\infty$ and applications to  signal processing.
{\em Problems of Information Transmission}, iss. 2. pp. 113 - 126.  
 
 \bibitem{D24p} Dokuchaev, N. (2024).  On  predicting  for  \bl{ non-decaying  unbounded } continuous time  signals. 	arXiv:2405.05566 
\comm{  \bibitem{D24b}  Dokuchaev, N. (2024). 
 Sampling Theorem and interpolation formula for non-vanishing
signals. arXiv:submit/5600348}
\bibitem{Eldar}
Eldar, Y.C.,  and T. Michaeli, T. (2009). Beyond bandlimited sampling. {\em IEEE Signal Processing Magazine}, vol. 26, no. 3, pp. 48--68. 


\bibitem{FeiMD}
Feichtinger, H.G.  (2024). Sampling via the Banach Gelfand Triple. In: Stephen D. Casey, Maurice
Dodson, Paulo J.S.G. Ferreira, and Ahmed Zayed, editors, Sampling, Approximation, and
Signal Analysis Harmonic Analysis in the Spirit of J. Rowland Higgins, Appl. Num. Harm.
Anal., pp. 211-242. Cham: Springer International Publishing.
 
\bibitem{GenH}
 Haddad, A. H., Thomas, J. B.,  Yao, K. (1967). 
 General Methods for the Derivation of Sampling Theorems. IEEE Transactions on Information Theory, 13(2), 227-230. 
\bibitem{Han}
Han, Y.. Liu, B., and Q. Zhang, Q. (2022).
 A Sampling Theory for Non-Decaying Signals in Mixed Lebesgue Spaces $L_{p,q} (R^{d+1})$, Applicable Analysis 101, no.,173–191.
24-25.

\bibitem{Kahane} Kahane, J.-P. (1970).
Séries de Fourier absolument convergentes. Springer Berlin, Heidelberg.
\bibitem{Katz} Katznelson, Y. (2004). An Introduction to Harmonic Analysis. 3rd Edition, Cambridge University Press, Cambridge. 
\bibitem{KPT}
 Kircheis, M., Potts, D., Tasche, M. (2024)  On numerical realizations of Shannon’s sampling theorem. 
 Sampl. Theory Signal Process. Data Anal. 22, 13. 
 \bibitem{Li}
Li, W., Wang, J. (2025). Random Sampling and Reconstruction of Non-Decaying Signals 
 From Weighted Multiply Generated Shift-Invariant Spaces. Math. Meth. Appl. Sci. 
 
\bibitem{NU} Nguyen, H.Q.,  Unser, M. (2017).
A sampling theory for non-decaying signals.
Applied and Computational Harmonic Analysis 43(1), pp. 76-93. 
 

\bibitem{RadExp}
Radchenko, D., and M. Viazovska, M. (2019). Fourier interpolation on the real line. Publ. Math. Inst.
Hautes ´ Etudes Sci. 129, 51–81.
\bibitem{RamExp}
 Ramos, J.P. G.,  Sousa M. (2023).
Perturbed interpolation formulae and applications.
 ANALYSIS \& PDE, Vol. 16, No. 10, pp. 2327–2384
\bibitem{Rei}
Reiter, H. and Stegeman, J.D. (2000). Classical Harmonic Analysis and Locally Compact Groups.
2nd ed. Clarendon Press, Oxford.

 \bibitem{Tr} Tropp,  J.A., Laska, J.N., Duarte, M.F., Romberg, J.K.,  Baraniuk, R.G.  (2009).
Beyond Nyquist: Efficient sampling of sparse bandlimited signals. IEEE Trans. Information Theory 56 (1), 520-544. 

\bibitem{VaaExp}
Vaaler, J. D.  (1985). Some extremal functions in Fourier analysis, Bull. Amer.Math. Soc. 12,
183–215.

\bibitem{V01} 
Vaidyanathan, P. P. (2001). Generalizations of the Sampling Theorem:
Seven decades after Nyquist. {\em
IEEE Transactions on circuits and systems - I: fundamental theory and applications}, v. 48, no. 9,
1094 --1109.

\bibitem{Z} Zayed, A.I. (1993). Advances in Shanonn’s sampling theory, CRC Press, Boca Raton, FL.

\end{thebibliography}
 \end{document}